# A multi-strategy optimizer for arbitrary generic functions in multidimensional space

Glauco Masotti [(*)]

**Abstract**

An algorithm capable of finding a likely global optimum (minimum) and a set of sub-optimal points for arbitrary generic functions of several variables is presented. The algorithm is designed to deal even with functions of complex behavior, irregular and noisy, with steep variations and exhibiting a lot of local sub-optimal points. The complications of having to deal with a finite domain, as this is usually the case, are taken into account. The method is composed of a number of cascaded stages, each employing a different strategy to improve over the results of the previous stage. Many ideas and concepts employed in known methods are re-elaborated in a coherent scheme, plus several new ideas are introduced. Line minimization plays an important role in most stages, for this purpose a new and powerful algorithm for line minimization is used as well.

## 1    Introduction

Optimization problems are commonly schematized mathematically as the problem of finding the global minimum of a function of one or more variables [19], usually referred also as parameters of the function.
For regular, well defined, analytical functions some well studied methods exist for finding extremal points [1], in these cases derivatives may be used for better computational performances.
Unfortunately the functions with which we have to deal in practice are usually much more complex and cannot be expressed in analytical terms.
This is the case for most practical application cases, where the functions to be optimized express a "cost" of a solution, or otherwise a "merit", or "goodness" of behavior of the process at hand. Processes may be of any kind: technical or industrial, i.e. related to complex artificial systems, financial processes, e.g. trading systems, or natural and biological processes. If a process can be modeled mathematically, and a function of cost or merit can be defined, then optimization of the process should be possible.
In all these cases the optimum point of the function represent in substance the best way of doing things, that is consume less energy, minimizing costs, maximizing gains, assure best probability of success to a therapy, and so on. It is thus evident the importance of having good optimization methods.
Nevertheless this is still an open field of research. No superior method for the optimization of generic functions in multidimensional space is recognized, all proposed methods have their pros and cons and are mostly based on heuristics. In fact theoretical results valid for the general case can hardly be defined or stated. Profs of convergence for the various proposed methods usually do not exist.
What can be done in practice is to find methods which increase the probability of finding the optimal point (or points) without increasing too much the cost of the method. Assuming that evaluation of the function is the most computationally intensive task, the cost of a method can be expressed with the number of evaluations (i.e. of points in the parameter space where the function is evaluated). The efficiency of a method can thus be expressed as $P(opt)/nf$, where $P(opt)$ is the probability of finding the optimal point $opt$, and $nf$ is the number of evaluations of the function. Usually we have to choose a trade-off between these divergent requirements, maximizing $P(opt)$ and minimizing $nf$. Most efficient methods should give the best results for the same cost, however this is not always true. In fact $P(opt)$ has an unknown form, it depends on the method used and is usually a non linear function of $nf$. For some methods it can "saturate" at a certain point, i.e. increasing $nf$ does not increase $P(opt)$ enough or anymore. Thus chasing efficiency may impede finding $opt$ in critical cases. In other word slow, costly methods may guarantee a result in the long run, while supposedly more efficient methods may give a quick result in simple cases, but not a good result in more complex ones.
An important property of a good optimization method should thus be not to exhibit such a saturation effect, i.e. increasing $nf$ should always increase $P(opt)$, this means that the optimizer should not be "trapped" in a sub-optimal solution from which it cannot escape. Some extra cost for simple cases can be sacrificed in order to guarantee this property in more critical ones.
However every method that can be devised must rely on some sort of regularity of the functions that have to optimized. At one extreme optimal points for very regular functions can be found exactly and easily with few steps of specialized algorithms, e.g. if the function is a paraboloid with circular cross-section just one step can be necessary! At the other extreme, if a function is pure white noise then probably any algorithm could hardly compete with random walk or exhaustive search.

For every method used the attainable results depend largely on the context and the particular problem at hand. Given the vastity of the literature on this subject a comprehensive review of known methods is beyond the scope of this paper, however references to known works which are relevant to our approach will be given and discussed throughout the paper.

## 2  Our approach: the multi-strategy optimizer.

The entire method which has been implemented in our optimizer is composed of a number of cascaded stages, each employing a different strategy to improve over the results of the previous stages. The different characteristics of the methods employed in each stage, as well as other provisions which have been adopted, should guarantee convergence even in most critical cases without compromising too much efficiency in simpler ones. This strategy is in general effective in avoiding that the optimizer is "trapped" in a sub-optimal solution, in fact the number of options and of possible paths of execution is very large and where a method fails another may succeed, so that chances of finding the optimal solution are greatly enhanced.

We used as test cases the function of merit of an experimental trading system with various sets of data. Any consideration on the particular application case is out of scope for this paper. The system behavior, and consequently performance, depends in our example case on 11 parameters, the solution space is thus an $11^{th}$ dimensional space, but the optimizer has been tested also with systems depending on up to 18 parameters. The objective function is a function of merit of the system, which takes into account the total gain, duration ad entity of drawdowns, distribution of gains over time, etc. The functions resulting from different datasets usually have a clear global minimum and a number of suboptimal solutions. However they are quite noisy and exhibit a lot of local minima which disturb the search for the really good regions. Nevertheless, as it will be shown, we obtained good results.

A precise comparison of the performances of this algorithm with respect to other known methods has not been attempted, however this is the final result of an evolutionary process which started some years ago, during which several algorithms have been tested.

Most of this work has been discarded, because the results were not satisfying, some has been kept and constantly improved, leading to the algorithm which is presented here.

We started implementing and testing well known classic methods [1], and then elaborating them with modifications and hybridization with other techniques and introducing new ideas as well.

In particular a considerable amount of time has been spent with an algorithm derived from the downhill simplex method [1]. At the end of this path of work the resulting algorithm was a lot different from the original formulation by Nelder and Mead [12],[13], in practice only the fact of using a simplex was saved, but at a certain point also this appeared as a limitation rather than an advantage, so that this "path of evolution" was terminated. However not all of this work went wasted, because it served to test many ideas which are still present in other forms in the final version of this algorithm. The algorithm presented here gave the best results in our tests, being able of finding the global minimum in all the test cases in a reasonable amount of time (see ch. 4), whereas other methods failed, even increasing the computation times considerably.

These results are referred to the serial version of the algorithm. A parallel version, which exploit parallelism at a "macro" level is under development, this should ensure computation times inversely proportional to the number of CPU's used.

We do not see reasons why our method should not give good results also in other contexts, however the algorithm behavior depends on many parameters, which should perhaps be tailored for each particular application for best results. Given the complexity of our approach, we think that our method is best suited for particularly "tough" cases, where simpler methods fail or do not guarantee of finding the optimal solution with a sufficient degree of confidence. Systems with a "noisy" objective function, exhibiting lots of local minima, and/or with the optimum region confined in a small portion of the parameter space (something which becomes more likely increasing the number of dimensions), are good candidates.

### 2.1  Basic concepts

*Normalized space*
In our method we assume that the domain of each variable (or otherwise said parameter) is a finite interval of the real numbers, therefore the entire domain of the objective function is an hyper-block, which can be mapped in the unit hyper-cube, which is assumed as the normalized domain of the objective function. We call this "normalization" of the parameter space. Each normalized parameter of the function is thus defined in *[0, 1]*.

*Swarm/stack*
The algorithm maintains an ordered set of the best distinct points found during execution. This set is called "stack" or alternatively "swarm", in fact our methods borrows some ideas from Particle Swarm Optimization (PSO) [14], [15]. In particular the set of best points found can be seen as a swarm, but this swarm is maintained ordered in a stack of finite length *nStack*. When the stack is full the worst point is eliminated, to make room for the last point found, if this is better. We experimented with *nStack* in the range *[40, 200]*, a final value of *120* has been used for our application case.
Not all points computed by the algorithm are placed in stack, e.g. when line minimization is performed [20] a number of points lying on the line are computed, but only for the best point found the placement in stack is attempted.

*Equivalent points*
The algorithm during execution computes many points in parameter space where the function is evaluated. Points which are close enough to each other in parameter space should be considered as "equivalent". We do not want the swarm to concentrate too soon in a particular zone, at least not before having explored enough the entire space, therefore only the best instances of these equivalent points are kept in stack. This constraint is relaxed as the algorithm progresses, which means that the "radius of equivalence" *rEq*, i.e. the threshold distance between two points below which the two points are considered as equivalent, is progressively lowered. This means that at the beginning the swarm is formed by well distinct points, while towards the end of the algorithm points are allowed to concentrate around optimal solutions, this allows for extensive exploration of space at the beginning, while towards the end attention is concentrated in the good spots, where the global optimum is likely to be located.

*Norm to be used*
Recognizing equivalent points requires the computation of a distance between the two points, so we have to decide which norm to use for computing this distance. The choice in favor of Euclidean distance $D_2$ is the most obvious, but may not be the best one. We experimented with various norms and decided to adopt the $D_1$ norm (Manhattan distance). This choice is to favor the presence in stack of points which differ in a larger number of parameters rather than in just one or few parameters. In fact at the other extreme we have the $D_\infty$ norm, where only the maximum difference in the parameters counts. For instance, using this norm in $R^3$ the points (0,0,1) and (1,1,1) have the same distance, equal to 1, from (0,0,0), therefore taking *rEq*=2 they will be declared all equivalent, while using the $D_1$ norm the distances evaluate respectively to 1 and 3, thus the points (0,0,0) and (0,0,1), which differ in only one coordinate, will be declared equivalent, while the points (0,0,0) and (1,1,1), which differ in all the 3 coordinates by the same amount, will be declared as distinct. The $D_2$ norm situates halfway, so it is clear that the desired effect is maximized using the $D_1$ norm.

*Temperature*
We borrow from simulated-annealing methods [1], the concept of "temperature" of the optimizer, but the way this concept has been implemented here is somewhat original.
At first the temperature *T* of the system is high, specifically we start with *T=1*. At high temperature points of the swarm usually make large moves at each step, i.e. the optimizer uses larger steps, that is it applies larger perturbations to the points already found to produce the candidate points for the next step. This means that we favor the exploration of space, in search of good candidates for the optimum, rather than concentrating in "refining" already found solutions (considering that a point in stack can be viewed as a solution of the optimization problem).
When the optimizer proceeds in its execution, the temperature *T* is gradually reduced in a sequence of discrete steps. With *T* approaching 0 the optimizer takes smaller and smaller moves. This means than towards the end of execution of the algorithm we do not waste time anymore in exploration, we focus on already found solutions instead, trying to refine them, pushing them closer and closer to the real local optimum in the vicinity.

*Line minimization*
Line minimization, i.e. one-dimensional optimization [1] is a fundamental building block of our optimizer as it is used extensively in many parts of it. We called *linmin* the procedure performing line minimization. We developed for this purpose an original algorithm which improves over Brent and similar methods [1],[16] in several aspects. It manages multiple local minima, takes into account the complications of having to deal with a finite domain, rather than an unlimited one, and has a slightly faster convergence in most cases. A detailed description of this algorithm is given in a dedicated paper published in concomitance with this one [20].
In general the line along which we have to perform minimization is oblique in the domain, i.e. it is not parallel to a coordinate axis, so that first of all we have to determine the domain *[a, b]* of the parameter *t* which specifies a point *P* along the line, given that for *t=0* we have the starting point *P0*. The parameter *t* corresponds to *|P – P0|*, the length of the segment joining *P* with *P0*. The limits of the domain *[a, b]* are calculated with a bisection algorithm as the smallest and largest numbers for which the point P is contained in the unit hyper-cube, which is, as we said, the normalized domain of the objective function.

*Distributions of probability*
The optimization algorithm is essentially not deterministic, its functioning depends in many of its parts on random choices, the nature of the problem is such that finding an optimal solution is not guaranteed, we can only make it very likely. In order to maximize the probability of this event it is essential to use suitable distributions of probability for the various stochastic processes which constitute the algorithm.

We introduce here for the first time some distributions of probability which are used in various parts of the algorithm.

*Generation of random numbers*

For the reasons exposed above it is essential that the generation of pseudo-random numbers be of very good quality! The generation of random numbers of any distribution is based on the generation of uniformly distributed random numbers, therefore it is essential that this basic function be performed at best. The rand() routine of the standard C library is not up to the task, first of all, but not only, because of the short period. Several alternative procedures for the generation of uniformly distributed random numbers, known for having a very large period and excellent statistical properties [8],[9],[10],[21],[22],[23], were considered. The differences are subtle and difficult to evaluate [23]. Therefore we proceeded to test the obtainable results. The generated numbers were mapped in the unit interval, and the couples made up with two consecutive numbers were mapped in the unit square. The uniformity of the generated distribution of points in several tests were then subjectively evaluated. We ended up choosing David Jones's JKISS [10] as the basic random number generator for our algorithm.

*Exploiting multiple trials*

The solutions found by the algorithm obviously depend on the initial guess, and if the initial guesses are the same, depend on the sequence of moves performed, which in their turn are basically determined by (pseudo)random processes. In the presence of a sufficiently complex function, with many local optimal points (that is to say many potential optimal solutions), if we repeat the entire optimization process (trial) with different random sequences, it is very likely that we obtain different results, i.e we should find final stacks with different solutions. If the optimal region is difficult to find, we may thus increase the probability of finding it by increasing the number of trials. This is in fact what we do, we perform several trials of optimization and we merge the results.

We use typically ten optimization trials in our latest implementation, which are performed at the same temperature *T*. After the completion of all sets of trials, all solutions from all available stacks are merged in one stack, by taking the first *nStack* best solutions. If *T>0*, these solutions are then used as guess points for the next cycle of optimization at a lower temperature.

It is relevant to point out that the execution of each trial is independent from the others, thus they can be executed in parallel! We think that the concurrent execution of the trials is the most natural and probably the best way for parallelizing the algorithm.

*Stack score: evaluating compliance of the swarm to the optimization task.*

Some measures have been defined and used to evaluate the compliance of the swarm to the optimization task and thus ascertain a good behavior of the algorithm in all of its phases.

We should avoid that the swarm prematurely focus on a sub-optimal solution, i.e. that all points in stack concentrate too early in one region or in a few regions. The swarm, particularly at the beginning, when *T* is high, should be dispersed all over the parameter space, in order to accomplish a satisfying exploration of all regions.

As a measure of dispersion we take:

$$statDist = ((\sum_{ij} d_{ij}^{1/2})/(n*(n-1)))^2 \qquad (1)$$

that is the square of the average of the square roots of all distances between the *n*(n-1)* couples of points of the swarm, i.e. a generalized mean with *p=1/2* of these distances (we called it *"square root mean"* as opposed to the more popular *"root mean square"*).

This measure assigns a penalty to multimodal distribution of points, while a quadratic mean does not. This can be easily verified in some simple 2D cases, as in the example shown in Fig. 1.

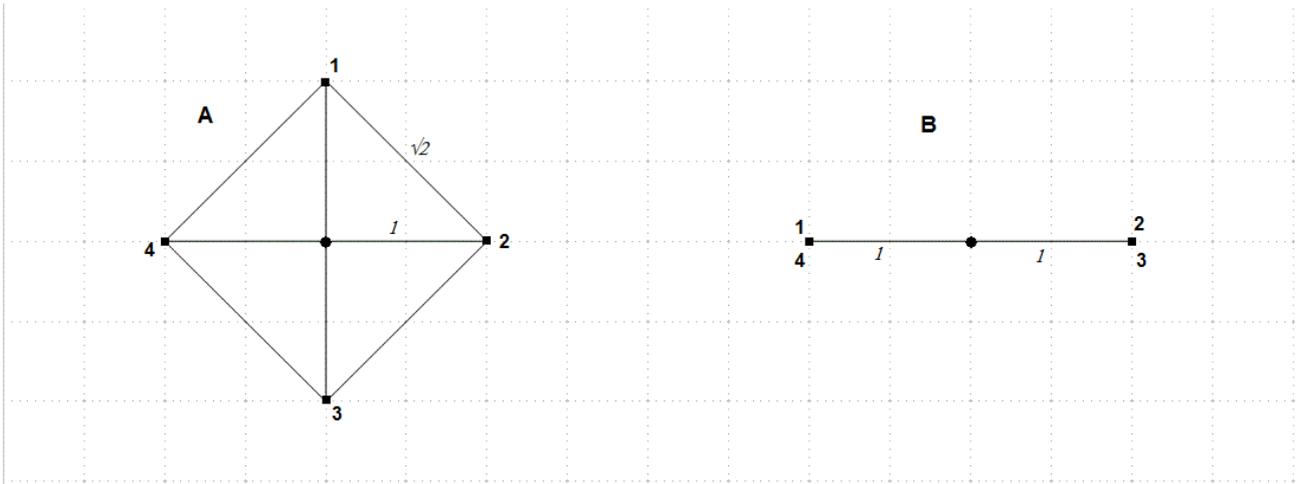

**Figure 1.** A simple example for testing different measures of dispersion.

Here we have four points placed at distance 1 from a center, in case A the points are placed at the corners of a square, while in case B they all lie on a horizontal line, points 1 and 4, as well as 2 and 3, are overlapped, thus they represent a sort of multimodal distribution, because the points collects in two loci in 2D space, instead in case A the points are more dispersed, thus we need to use a measure which penalizes case B. If we use as a measure the root mean square of distances, we have $((4*(\sqrt{2})^2+2*2^2)/6)^{1/2} = (8/3)^{1/2}$ for case A, and for case B: $((4*2^2+2*0)/6)^{1/2} = (8/3)^{1/2}$ ! Thus the measure does not distinguish between the two cases. It can be verified instead that using a simple arithmetic mean of distances case B is slightly penalized, and is severely penalized if we use the "square root mean" of distances.

More detailed statistics for each parameter are also collected. We should avoid the premature concentration of the swarm around particular values of the parameters. A sensitive measure, which degenerates to 0 when all points in stack have the same value for one or more parameters is the harmonic mean of the standard deviations of each parameter. Let us call this measure *statParams*. A cumulative measure of dispersion is given by:

$$dispScore = ((statDist)^{1/2} + statParams^{1/2})/2)^2 \qquad (2)$$

Using a square root mean also here makes the measure more sensitive to degeneration in just one of its components. A comprehensive measure of "fitness" of the swarm is provided by computing the exponential moving average of the values of the objective function (treating them as a time series, starting from the worst point), let us call this *fmtScore*. Finally an overall score of the swarm/stack is evaluated by putting everything together:

$$stackScore = fmtScore*(1+ dispScore) \qquad (3)$$

While the algorithm proceeds and temperature decreases, it is normal that the dispersion becomes smaller, because the swarm concentrates on good regions, while the fitness should monotonically increase. However, for a well working algorithm, the overall score should always increase. If this is not the case there should be something which is going wrong, probably some of the parameters or options controlling the behavior of the optimizer (as we will see) should be adjusted to fit the problem at hand.

## 2.2   First stage: swarm population and random search.

At the very beginning we should make the assumption that little is known of the behavior of the objective function to be optimized, and that we do not have guesses for the locations of the optimal points. Therefore we take some default points to start with. A simple but effective choice is to take three points on the principal diagonal of normalized space corresponding to the normalized coordinates ¼, ½, ¾. Fig. 2 shows the projection of these points over some planes of the 11th dimensional space of the test application problem.
We used these projections for debugging purposes. The unit vectors (versors) parallel to coordinate axes, or the linear combinations of these which have been used to define the orthogonal versors which specify the projection planes, are indicated.
To understand the graphics of Fig. 2 please take into account that planes which contain faces of the unit hypercube are represented without distortion, it is so that the projection of the sides of the hypercube result in a square, and they are displayed so. Planes which are oblique with respect to faces of the unit hypercube intersect it in a rectangle, these planes have their projection rescaled in the longer dimension, in order to be represented as a square as well. It turns out doing so, that the projections of the unit sphere are ellipses (displayed in gray), while the ellipses inscribed in those rectangles are instead projected as circles (displayed in black).

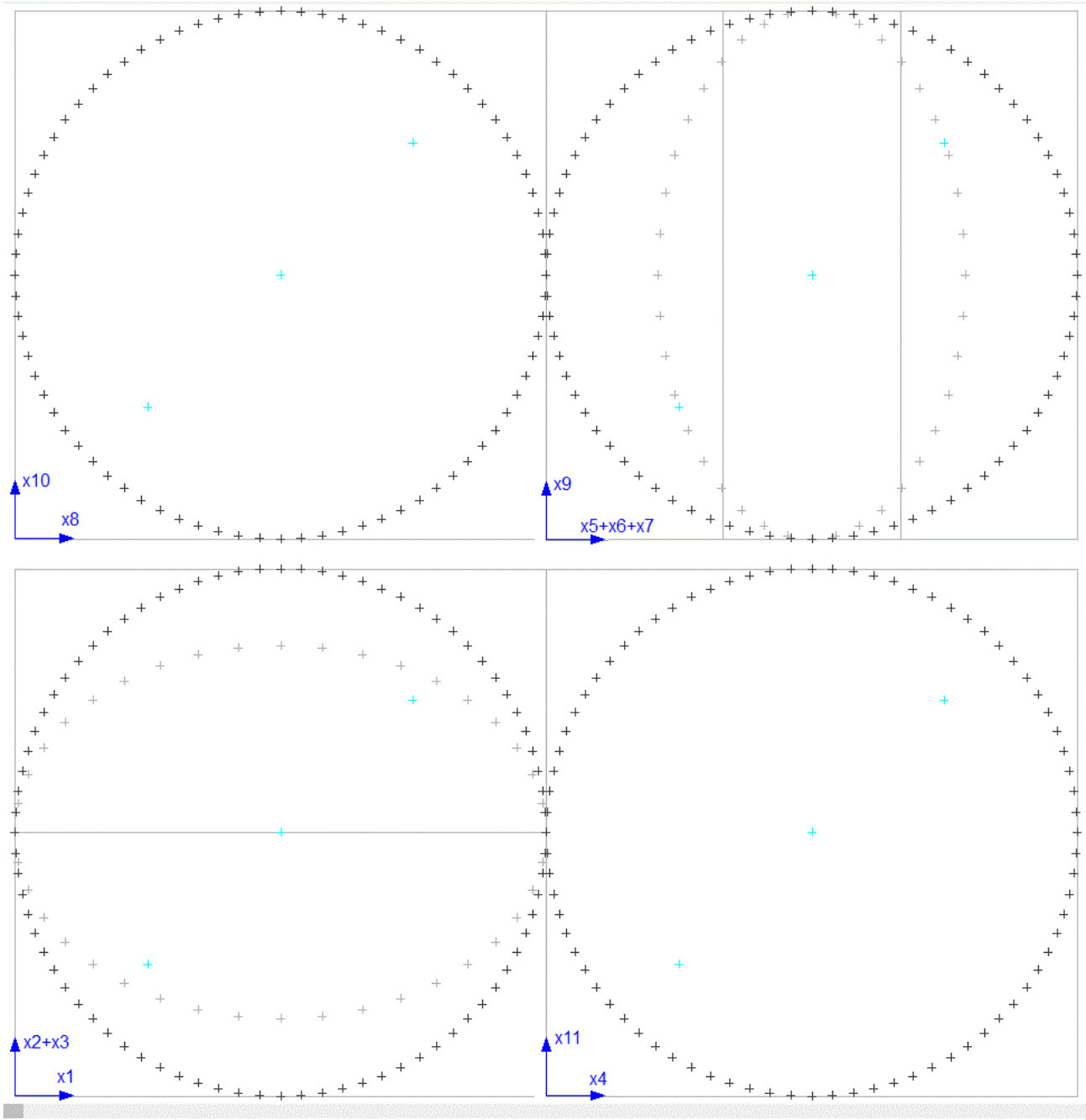

**Figure 2.** Projection planes used for debugging purposes and the versors which identify them.
The projections of the three default initial guess points along the principal diagonal
are also shown in cyan.

The purpose of the first stage of the optimizer is essentially to explore the parameter space. Thus random steps in some specified random directions are taken from the three initial points, in this way we arrive in three additional points.
The basic step is to add these points to the swarm, i.e. to the set of the best points found so far, then the procedure of generating new points originating from existing ones moving along a random direction with a random step is reiterated. The new points are added to the swarm/stack. It is so that the swarm is populated, it doubles at each iteration and reaches its maximum population quickly.
From this point on the stack is managed with a "better in worst out" policy.

This is the basic procedure of the first stage, that we called *"swarm search"*. Let us examine it more in detail and let us point out some notable exceptions which apply.

1) First of all, what happens if one of the new points found is close (or even coincident) to one of the existing points? As we said above, only the best instance of equivalent points is kept in stack.

2) What happens if the new point which is generated is better than the original point? If a new point *P1* is better than the original one *P0*, i.e. *f(P1) < f(P0)*, it means that along the line passing through *P0* and *P1* there exists an optimum point *Pbest*, where *f(Pbest) <= f(P1)*, therefore it may be convenient to find this point and not just be happy with *P1*. In the worst case we have *Pbest = P1*, but in general we may have the opportunity to improve considerably the solution found. This is left as an option, but we have found that the extra cost of *linmin*, i.e. of finding this minimum, using the algorithm described in [20], is in practice worth it. In this case we place in stack *Pbest*, rather than *P1*.

3) How to choose each time a random direction for the next step? Various alternatives have been explored. One obvious choice is to generate every time a different random direction. A part its computational cost this does not appear to be the best choice. In fact, when the swarm is fully populated, the probability of generating this way colliding search directions, which may generate equivalent points is not negligible. Another possibility is to generate directions which originates from a single point, e.g. the center of gravity of the swarm, thus the swarm will expand or contract around the central point, but in this way the region around the center becomes "over-explored", while regions at the border are explored less.
A third alternative, which resulted in better performances, is to chose each time a unique random direction for the entire swarm. A random direction has thus to be computed only once at every iteration, and all the points of the swarm move along parallel lines. This reduces the probability of generating colliding points. While the search direction is common for all points, the random step which is taken is different for every point, negative steps are allowed too, so points of the swarms may also move in opposite directions.

4) How to chose each time the step for each point? Various probability distributions have been tested, we obtained the best results with what we called the *"notchTwinPeaks"* distribution. At first we used a *"twinPeaks"* distribution, later we added a notch at the origin.
These distributions are an original development and were specifically designed to meet the requirements of the algorithm. For *twinPeaks* the probability density as a function of *x* is defined by the expression:

$$twinPeaks(x,s,kt,q,ks) = (1-q)*(Gauss(x,-kt*s,s)+Gauss(x,kt*s,s))/2+q*Gauss(x,0,ks*s) \qquad (4)$$

Where *Gauss(x,m,s)*, is the probability density in *x* of a Gaussian distribution of mean value *m* and standard deviation *s*. The distribution is thus the sum of three Gaussians. Fig. 3 compares the resulting distributions for various values of *kt*, ranging from 0.9 to 1.3, while maintaining fixed the other parameters (in this case s=1, q=0.3, ks=2.5).

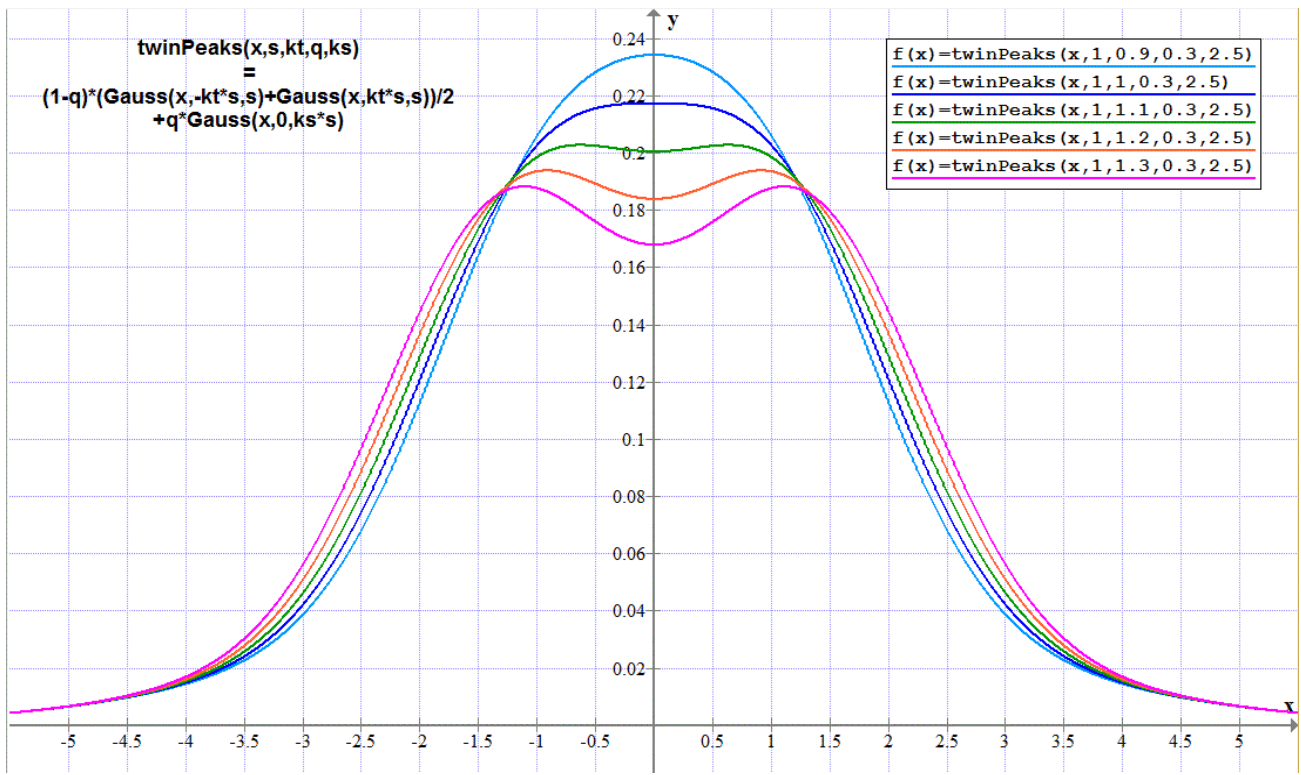

**Figure 3.** Effect of varying *kt* in the *twinPeaks* distribution.

For *kt=1* we have a flat top. Choosing *kt = 1.1*, the distribution is quasi-flat in a large neighborhood of 0, it has just two small "overshoots" on each sides, this considerably extends the extension of the plateau where the density function is quasi-maximal. Then the probability density decreases, but not as fast as a Gaussian does, finally the probability density for large values is still not negligible, which characterize *twinPeaks* as a "fat-tail" distribution.
The justification for taking the steps according to a distribution of this kind is as follows:
- The original point is one of the best point find so far, but it is not likely, especially at the beginning, that it is already the best point of this zone, it is more likely that it is located just in proximity of this (still unknown) local best point, therefore it makes sense to search for a better point not far from the original point. This accounts for the quasi-flat behavior around 0, with maximum probability density located here.
- However it seems reasonable to look sometimes somewhat further, because there may be better points that have not been discovered yet. This is the reason why the probability of taking larger steps decreases gradually, but it maintains a non negligible value also for quite large steps. This means that the algorithm is capable of taking sometimes what are called "Levy flights" [2], [6].

The standard deviation *s* is the main parameter that controls the dispersion of values around 0. Typically, in normalized space, we use *s = 0.05+0.35\*T*. Being the *twinPeaks* distribution, a linear combination of Gaussian distributions, we can either substitute *s* with the appropriate value *s* in (4), or use the distribution for *s*=1 and scale the values applying a factor *s*. When *T* decreases the values are scaled down, thus they concentrate more in a neighborhood of 0, this means that in general the optimizer takes smaller steps. Given these assumptions Fig. 4 shows how the *twinPeaks* distribution varies with respect to the temperature *T*.

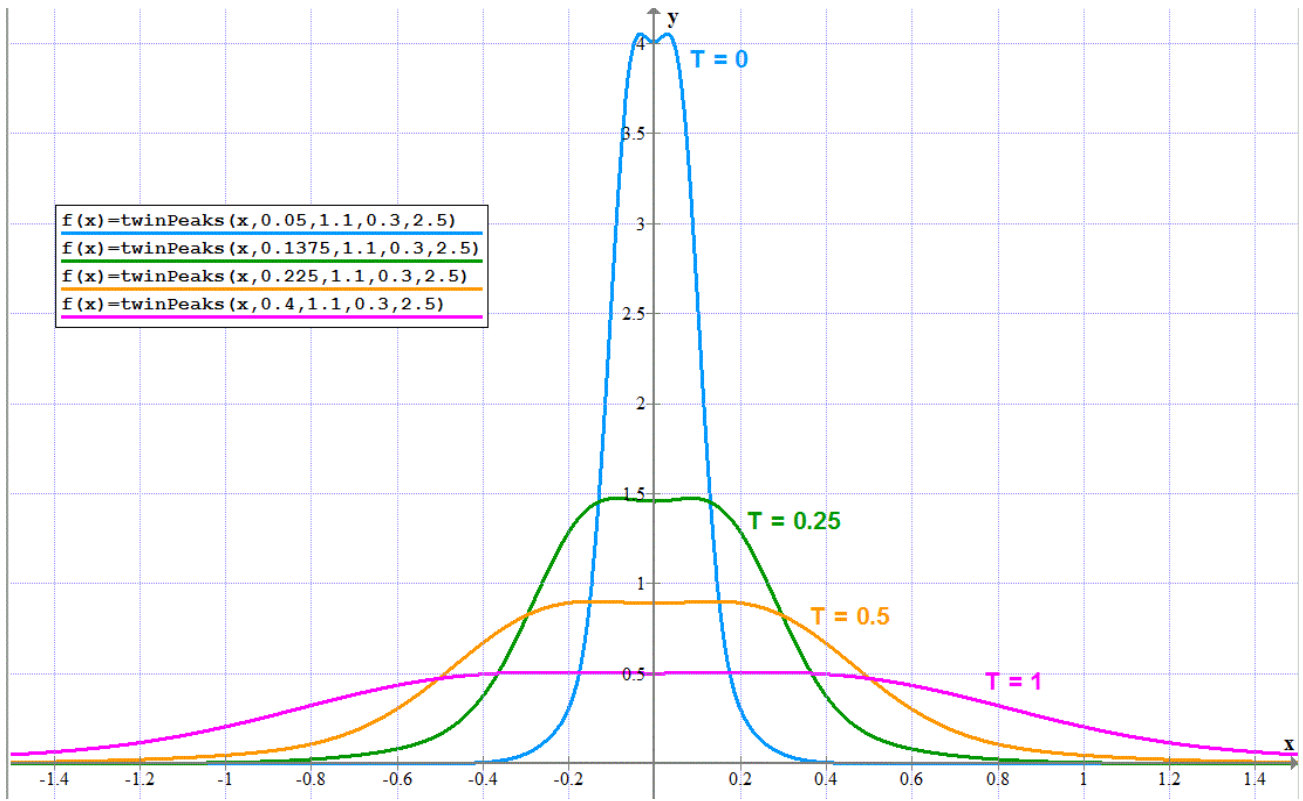

**Figure 4.** Variation of the *twinPeaks* distribution with *s(T)*.

It follows also that the generation of random values following the *twinPeaks* distribution can be made directly by generating numbers following each component Gaussian distribution in the appropriate proportions. However, in our case, we have to generate numbers in a limited domain, thus we must discard values outside the domain. If this is large, so that the probability of generating values outside the domain is small, the direct method can be used without problems, but if the domain is small, the direct method becomes inefficient. We developed for these cases a specialized procedure based on the accept/reject method. The final procedure which has been implemented uses one or the other of the two methods, choosing the most efficient for the case at hand.

The *twinPeaks* distribution however may often generate also very small steps, and this seems somehow a waste, especially at the beginning, because, unless the function is pathologically irregular, we should find similar values of the function stepping very close to the original point. The *twinPeaks* distribution has been therefore corrected by placing a narrow notch centered at 0. The final *notchTwinPeaks* distribution assumes the form:

$$notchTwinPeaks(x,s,kt,q,ks) = |x|^{1/16}*twinPeaks(x,s,kt,q,ks)/Area \qquad (5)$$

where *Area* is the integral of $|x|^{1/16}*twinPeaks(x,s,kt,q,ks)$ over the entire domain. The *notchTwinPeaks* distribution used for *T*=1 (which has *s*=0.4) is shown by a solid green line in Fig. 5.

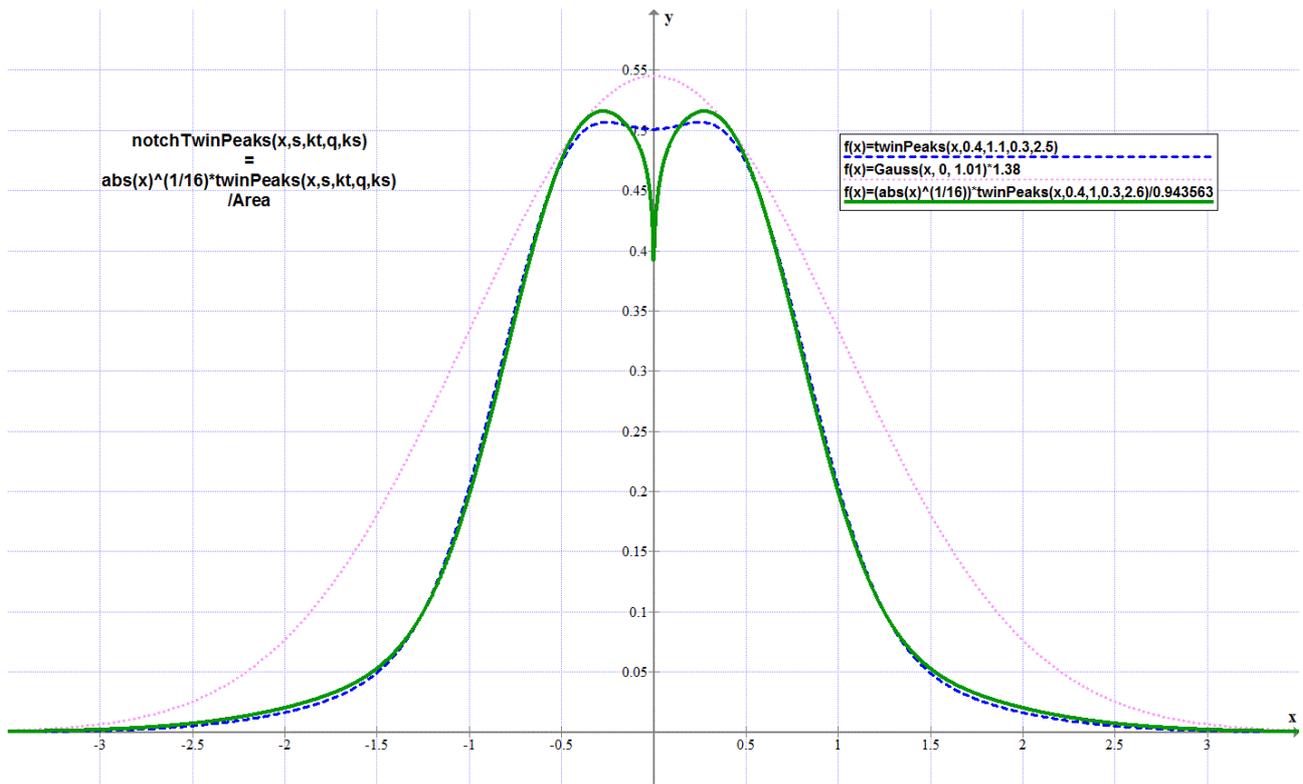

**Figure 5.** Matched *twinpeaks* and *notchTwinPeaks* distributions.

To compare the results with our previous experiments using the *twinPeaks* distribution, we looked for a close match between the two distributions by adjusting the other parameters. We see a close match with the *twinPeaks* distribution represented by a dotted blue line.

The distribution for *T=1* shown in Fig. 5 is taken as a reference. Results for smaller values of *T* are obtained by applying a scale factor to the numbers produced by this distribution. It turns out that, due to the correction factor $|x|^{1/16}$, which creates the desired notch at *x=0*, it is mandatory to use a specialized method for generating numbers following this distribution. An efficient accept/reject method has been implemented. Fig. 5 shows, with a dotted pink line, the "amplified" Gaussian function that can be used for this purpose.

Resuming our review of the first stage of our algorithm, Fig. 6 summarizes the effects of swarm search for our test case after its completion, as part of a trial with *T=1*. The projections of all points which belonged to the swarm, at least for a while, during the search is shown. There are also some groups of parallel yellow lines displayed, which correspond to the lines were *linmin* was performed.

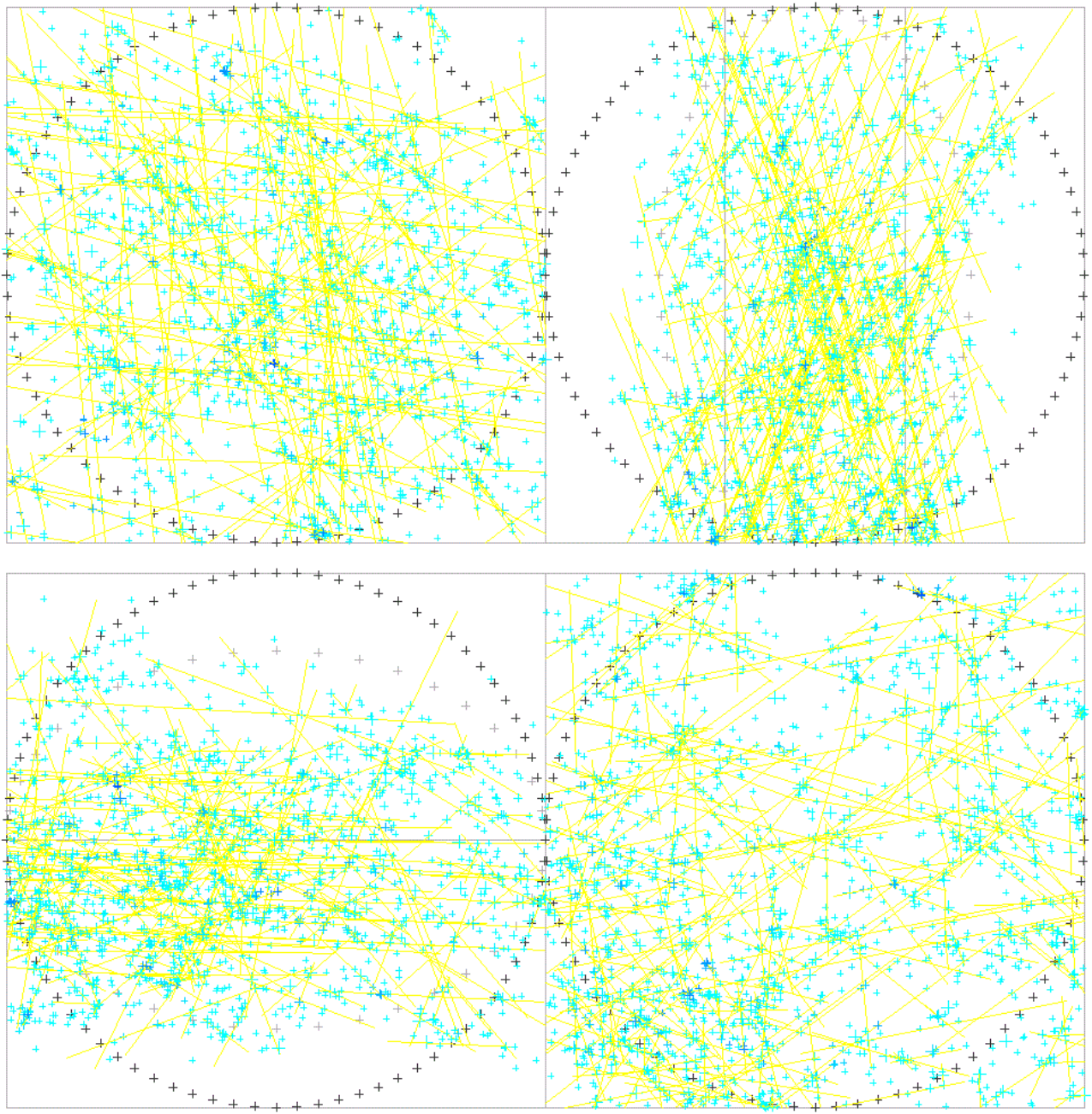

**Figure 6.** Stage one: swarm search

Points are depicted in pseudocolors, reflecting the associated value of the objective function, according to the color code shown in Fig. 7.

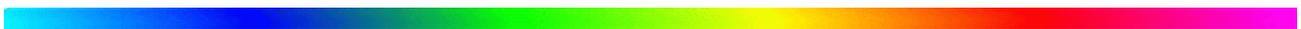

**Figure 7:** Pseudocolors used from higher (to the left) to lower (to the right) values of the function.

## 2.3    Second stage: genetics improvement

The second stage applies a genetic algorithm to the swarm to improve the solutions found. Genetics algorithm principles and mode of operation are well known [17], [3], however not enough theory supports the modeling of a problem and particularly what should be regarded as chromosomes and genes. These aspects are largely arbitrary. In literature most implementations tend to follow an analogy with the biological processes that inspired the method, thus the most used approach is to define as chromosomes collections of genes which are defined as sequence of bits at a very low level. Even if the problem is well modeled by a set of parameters which defines a particular solution of the optimization problem and these parameters are essentially floating point numbers, these are usually encoded in sequences of bits, which are split further to represent the genes (see for example [4]), regardless the fact that doing so genes may result in being meaningless, i.e. in not being directly related to any particular behavior of the system to be optimized. Genetic algorithms are powerful enough so that they work anyway, but we think that using "meaningful" genes would result in better convergence, and more control on the optimization process can be exerted as well. Following this reasoning we opted to keep the parameters themselves as genes, and the set of parameters defining a solution as a chromosome. This seems to us an obvious choice, nevertheless it does not seem to be so popular in literature, although some studies assess that using a floating point (FP) representation for genes, at least for the class of problems which are naturally modeled in this way, is advantageous [18].

Besides this, our genetic algorithm is peculiar in several other ways. The crossover of genes in our algorithm takes the form of recombination of parameters among two or more sets of parameters, i.e. a new set of parameters is produced by taking randomly each parameter from the generating sets. In a certain sense in this way we implement filiation from more than two parents, i.e. we do in one step what nature does in several generations. If the generated set inherits good characteristics from its parents and they do not interact in a negative way with each other, then the child set may improve over already found solutions.

Any parameter can also have a certain probability of mutation, which means that the algorithm is not limited to use only parameters already present in individuals of the swarm, but also has a capability of moving virtually anywhere in the solution space.

Temperature influences operations of our genetic algorithm in several ways. First of all we should note that, because the stack/swarm collects the best distinct solutions found so far, just the recombination of these individuals, treating them in an uniform way, causes a process of natural selection, because generated points which represent better solutions than existing ones will always find a place in stack, to make room, if necessary, the worst point will be ousted from stack. This process can be enhanced by favoring the recombination of the best points in stack. In our algorithm the formation of the set of points to be recombined is the result of several random choices. This architecture is largely arbitrary, it has been determined using good sense and heuristics, experimenting in the field.

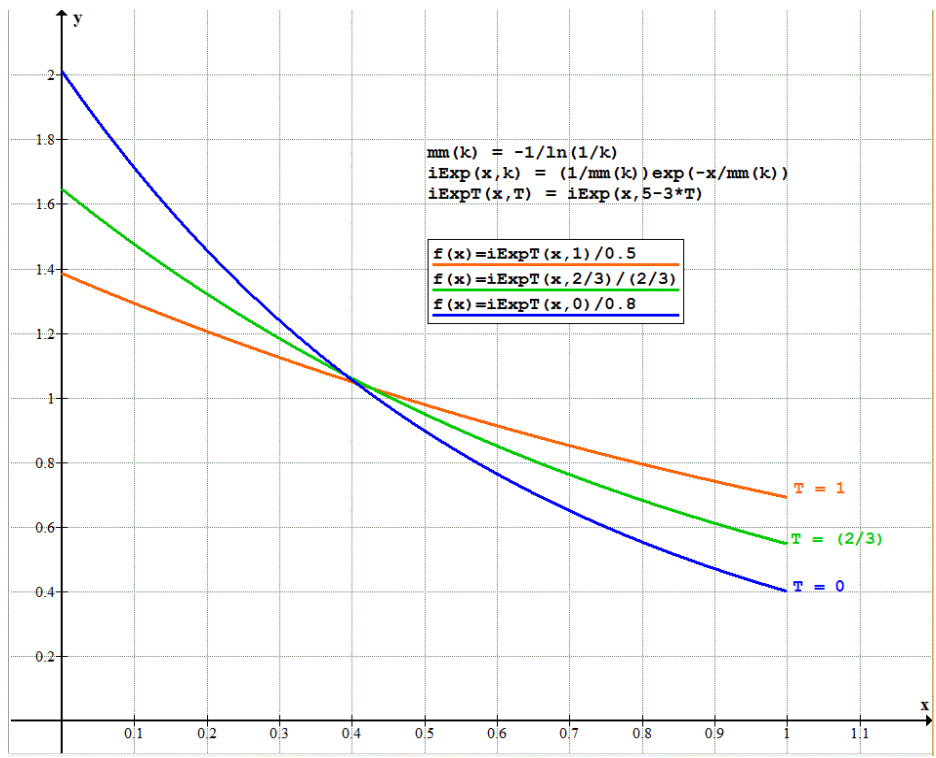

**Figure 8.** Bounded (truncated and rescaled) exponential distribution used in the recombination process

At first *nRecombine*, the number of points to be recombined, is calculated. A bounded exponential distribution is used to generate this number in the range *[2, nStack]*, where *nStack* in the number of points in the stack. The actual set of points to be recombined is formed by the first *nRecombine* points in stack, thus the best points are greatly favored and attention must be paid not favoring them too much, because this may cause to prematurely focus on a sub-optimal solution, where all points in stack will tend to concentrate. A chance, particularly at the beginning, when *T* is high, should be given to all points in stack, and these should be dispersed all over the parameter space. The parameter of the exponential distribution (the rate of the distribution, commonly called *lambda*) is a function of *T*, we assume *lambda* = 5-3*T*. Fig. 8 shows the resulting distributions over the interval *[0,1]* for *T=1, 2/3, 0*.
The unit interval is mapped in *[2, nStack]*. With these assumptions, when *T=1* the probability density for 2 is about two times than for *nStack*, while for *T=0* it is five times. Thus a small number of good "parents" is usually preferred, particularly at low temperature.
An analog exponential distribution, but mapped in *[1, nRecombine]*, is used to determine which value to pick up among those of the parents, for each of the parameters. Doing so, not only the best points have an higher probability of being parents, they also have an higher probability of transmitting a larger number of genes to their children.

After a parameter has been chosen, there is also a certain probability that it undergoes a mutation. We got best results with a probability of mutation of 0.25. At first we used a truncated *twinPeaks* probability distribution (but with *s(T)* resized to smaller values), centered at the original parameter value, to determine the mutated value. Fig. 9 shows an example of the distribution used for *T=1* in normalized coordinates. If we consider the total distribution of probability in *x*, with *x* in *[0, 1]* for a parameter of value *x0*, it assumes the form drawn with a blue line in Fig.8 (the figure shows the case *x0=0.3*), which is composed of a Dirac delta at *x=x0*, of area 0.75, together with a continuous truncated *twinPeaks* distribution of area 0.25.
This is a discontinuous distribution, but we can also imagine a continuous distribution approximating it. The infinite spike for *x=x0*, can be substituted by a narrow bell shaped distribution centered at *x=x0*, the result is a fat-tailed continuous distribution which assigns a probability of ~ 0.75 to values around x0, and ~ 0.25 for all other values. The purple line of Fig. 9 shows a distribution of this kind, which we called *fatTail3*, because is the composition of three Gaussians according to the formula:

$$fatTail3(x,c1,c2,c3,k1,k2,s) =$$
$$(c1/(c1+c2+c3))*Gauss(x,0,s)+(c2/(c1+c2+c3))Gauss(x,0,s*k1)+(c3/(c1+c2+c3))Gauss(x,0,s*k2) \quad (6)$$

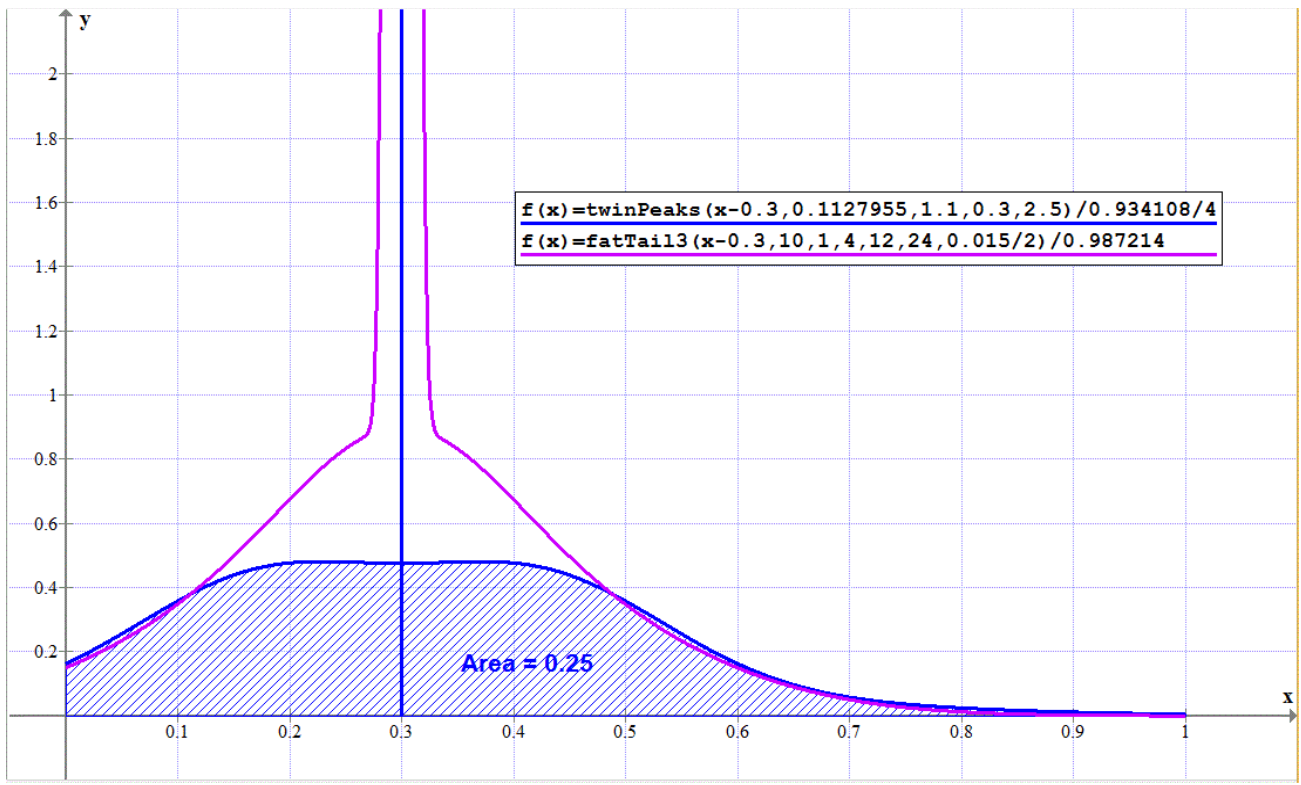

**Figure 9.** Two distributions tested for generating mutated values.

Choosing appropriately the parameters *c1, c2, c3* and *k1, k2,* the previous discontinuous distribution can be approximated and perhaps improved. In fact such a distribution would not commit strictly to the original value of the parameter, it will usually propose small variations of it. Doing so, in case the original value is close to the optimum, we have chances to improve the approximation. However also some large variations are proposed from time to time, doing so, in case the original parameter is far from the optimum, we may have some chances of stepping to the optimal region, but the probability of generating larger perturbations does not changes abruptly, there is an intermediate transition zone which can be useful when only moderate changes are necessary to reach the optimum. We tested various distributions of this kind, obtaining better results with respect to their discontinuous counterpart, *fatTail3* gave the best results in our tests and was adopted for the final version of the algorithm.

The generation of values with *fatTail3* distribution can be obtained by generating numbers following the component Gaussian distributions in the appropriate proportions, however a complication is given by the fact that we must generate in a limited domain. Here the same considerations made above for the *twinPeaks* distribution apply, thus, for this case too, we implemented an efficient procedure for generating Gaussian variables in a limited domain using the accept/reject method.

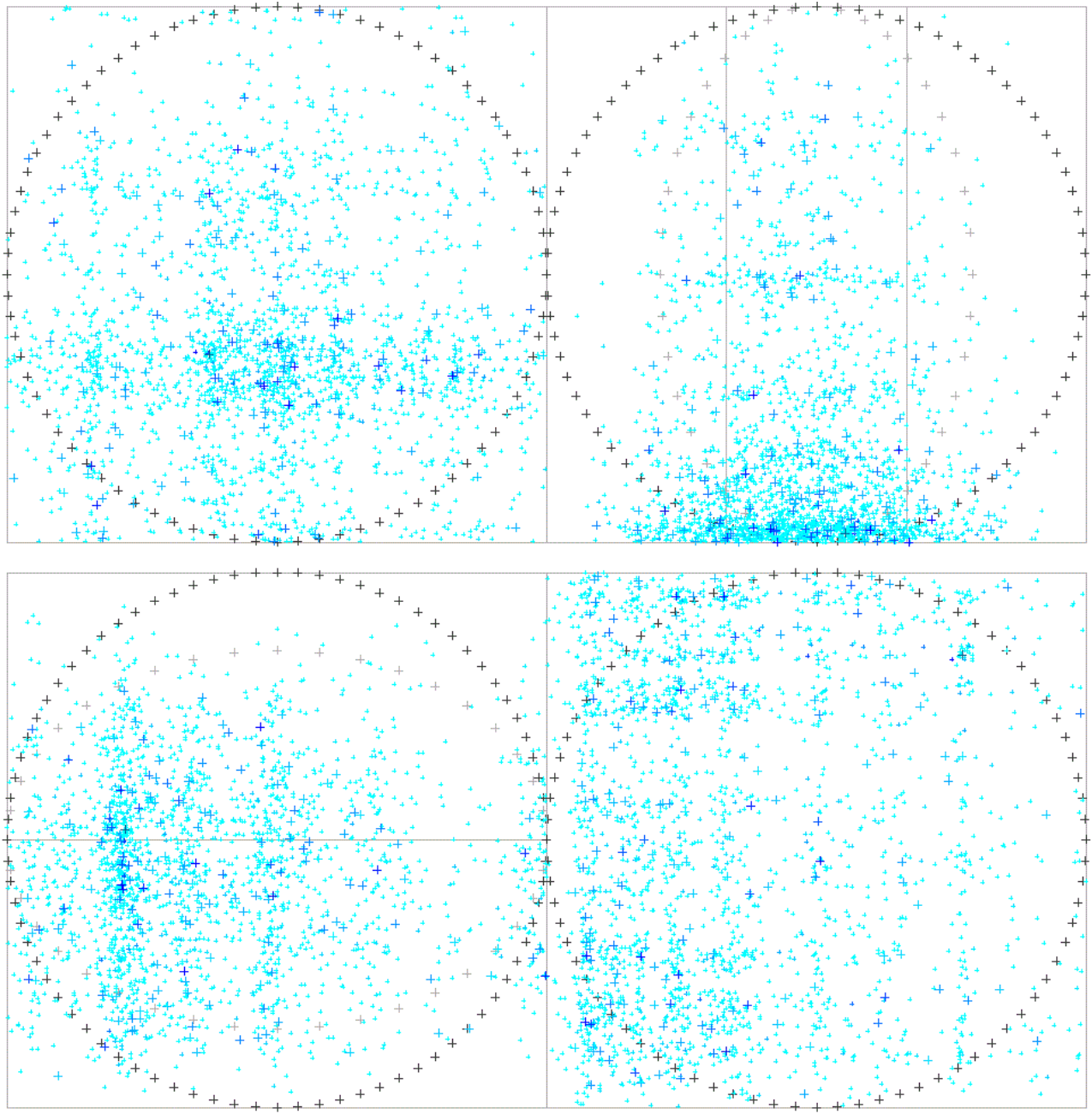

**Figure 10.** Points generated at an intermediate step by the second stage (genetic part of the algorithm) for *T*=0.75.

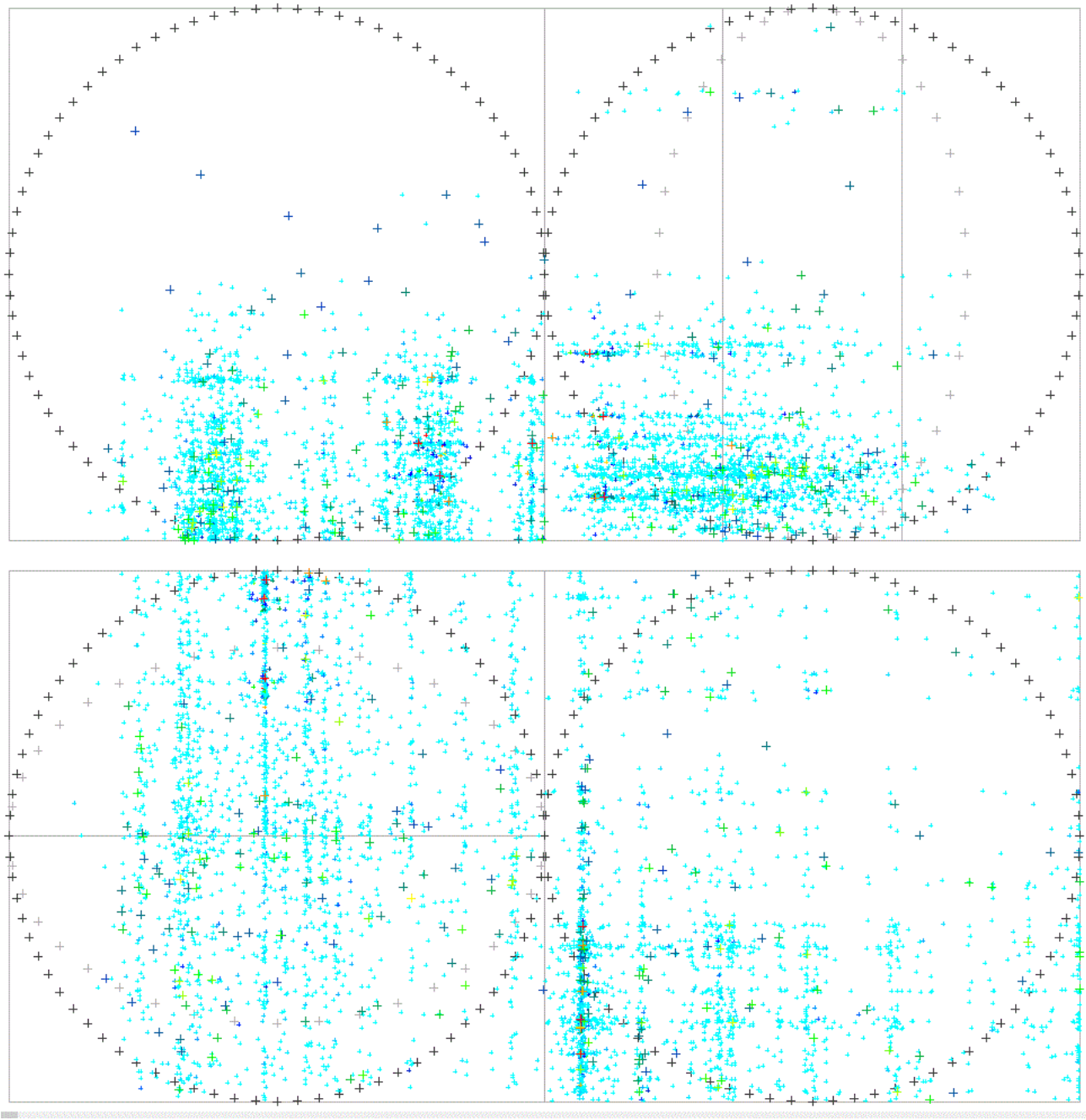

**Figure 11.** Points generated at an intermediate step by the second stage (genetic part of the algorithm) for *T*=0.

Fig. 10 and Fig. 11 show the projections on the test planes of all points generated by this stage during the optimization process at two intermediate steps, with *T*=0.75, and *T*=0. It can be seen how the generated points, as T→0, tend to distribute along strips having similar values for one or more parameters, which correspond to the current best choices for these parameters.

*2.4    Third stage: proximal search*

The third stage has been inspired by the Firefly algorithm [5], [6], but we think of having improved on it in several aspects, especially by exploiting line minimization, rather than simply taking random steps in the current search direction (this seems to be something not attempted in known approaches [7]). Because at this stage the fate of each point of the swarm is determined by the surrounding points, we have called this stage "proximal search".

The rationale of this part of the algorithm is as follows. We start from the worst point in stack and we search for a better point in the direction of one of the other points, we call this point the "attractor". Because all other points are better, there must necessarily be such a point in this direction, at worst it will coincide with the attractor itself, but doing *linmin* along the line joining the two points we should be able in general of finding a better point somewhere along the line. In this case the new point is added to the stack, otherwise, if no new point is found, we chose another attractor (until it is possible) and repeat the procedure.

The attractor point is chosen as the most "attractive" one, thus we evaluate a function of attractiveness for all points. In analogy to the Firefly algorithm, this function is proportional to the apparent luminosity of fireflies in a non perfectly transparent medium. Apparent luminosity is proportional to absolute luminosity and will decrease with distance for diffusion in space and also for adsorption of the medium. An isotropic diffusion of light in an n-dimensional space should cause a decrease of intensity as $d^{-n}$, where $d$ is the Euclidean distance from the source, but for the purpose of this algorithm we are not forced to follow this rule, in fact less severe attenuation could work better, severe attenuation may end up penalizing too much distant points, which may instead represent interesting search directions. A straightforward inverse square law has proven to work well (this has an analogy in spaces of dimension n >2 with the directional diffusion of light strictly in a plane, however by rotating it around a variable axis spanning the entire space is possible). An inverse quadratic exponential function of distance can simulate well adsorption of the media. Instead of using a constant attenuation coefficient, as described in literature [5], we use a coefficient which varies randomly, simulating variable conditions of visibility. Therefore, even in identical configurations of the swarm, the attractor point will vary, depending on visibility conditions. With low visibility, like when foggy, closer points are favored, instead when clear also distant but bright points may become attractors. This shrewdness enhances the number of search directions which are tried, and therefore the likelihood of finding the global minimum.

The attractiveness function should also depend on temperature, in fact when $T$ is high the algorithm should favor exploration of space, rather than local optimization, and vice-versa. We can do this by assuming in prevalence a "clear sky" when $T$ is high, and a "foggy whether" when $T$ is low. The final attractiveness function assumes thus the form:

$$attr(d) = a0*((1-T)*exp(-(d/D)^2) + T/(1+(d/D)^2)) \qquad (7)$$

where *a0* is the attractiveness at 0 distance (we can assume, for a point *p*, *a0=f(p)*), and $D$ is the characteristic distance, which is the inverse of the attenuation coefficient, the larger is $D$ the smaller is the attenuation with distance. As said, $D$ is a random variable, recomputed at each step of the algorithm. We used for this purpose the truncated Gamma distribution of Fig. 12.

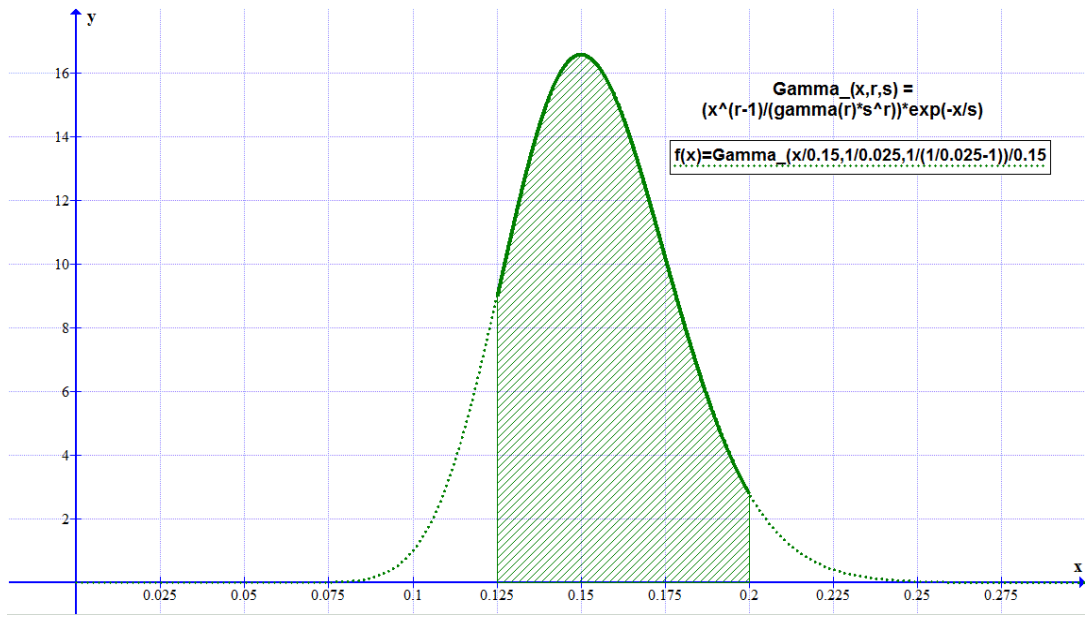

**Figure 12.** Distribution used to generate a random value for the characteristic distance.

The effect of this stage is illustrated in Fig. 13. It shows the projections on the test planes of all the points generated by this stage during the optimization process at an intermediate step with $T$=0.75. Also the traces of the lines representing the search directions were *linmin* was performed are shown. Most of these lines are shown in yellow, but some of them are drawn in magenta. These are the search directions which are attempted more than once. In fact it would be most likely unuseful to execute *linmin* two or more times on the same line, therefore each search direction is checked against the latest directions already tried and if the current direction is not new, *linmin* is skipped. It can be seen that a lot of directions are examined, but that there are some points that acts as the best attractors, so that most of the search directions are towards them.

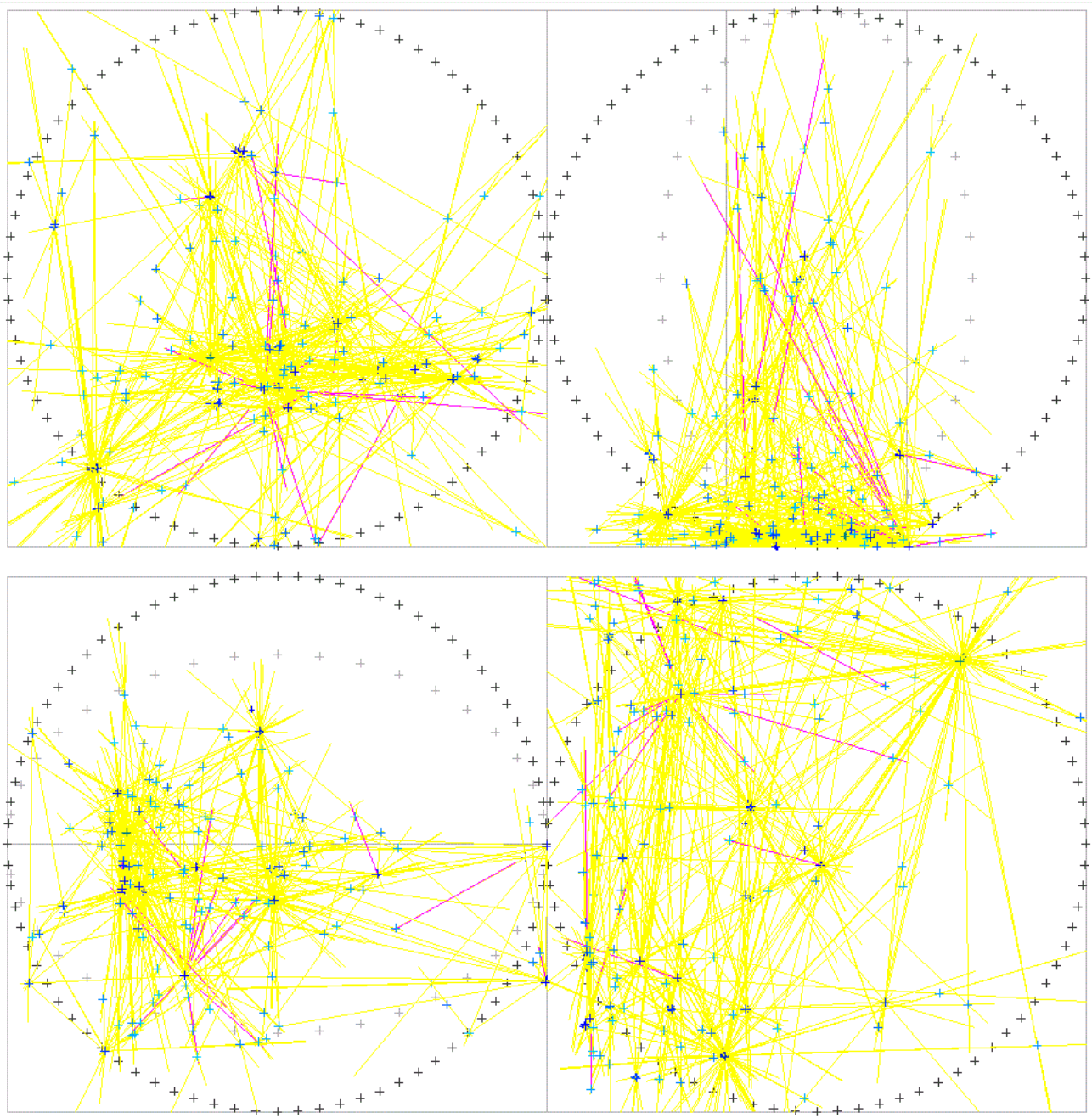

**Figure 13.** Projections of search lines and generated points during execution of the proximal search stage.

## 2.5  *Fourth stage: swarm along axes*

The fourth stage has been called "*swarm along axes*", because it does exactly this: instead of moving the swarm along random directions in parameter space, it moves the swarm only along the coordinate axes, but in a random sequence however. This means that at each step all points in stack will be varied in the value of just one parameter, while leaving the other ones unchanged.

Sequential search of a sub-optimal solution moving along coordinate axes is a simple but powerful idea, which has been used since early ages of multidimensional optimization, and has been proposed since then in many variants. One of the very first algorithm of this kind was proposed by Fermi and Metropolis [11]. Another variant which uses line minimization along each coordinate axes is described also in *Numerical Recipes* [1] (since the first edition), in chapter

10.
Although for some classes of functions the use of conjugate directions, like in Powell method [1], leads to a faster convergence than using coordinate directions, this is not true for generic and highly irregular functions, therefore we simply preferred to use the coordinate directions. After all these should represent privileged directions, otherwise we should not have identified a certain entity as a parameter of the system.

Our method applies line minimization along axes parallel to one of the coordinate directions, so that only one parameter is allowed to vary in this one-dimensional minimization process. Which parameter has to vary is determined according to a random sequence and not a predetermined one (e.g. the natural sequence in ascending order). This sequence determines the path which is followed in parameter space to arrive at the current minimum. In general the minimum found would be path dependent, thus using a random sequence should help in exploring space more effectively.

Known approaches deal only with a point at a time, i.e. starting from an initial point the current optimal point is moved at each step to a better position moving along the coordinate directions.

Applying this scheme to our swarm, this process would be repeated for every point of the swarm, taking one of them at turn as an initial point for the optimization process.

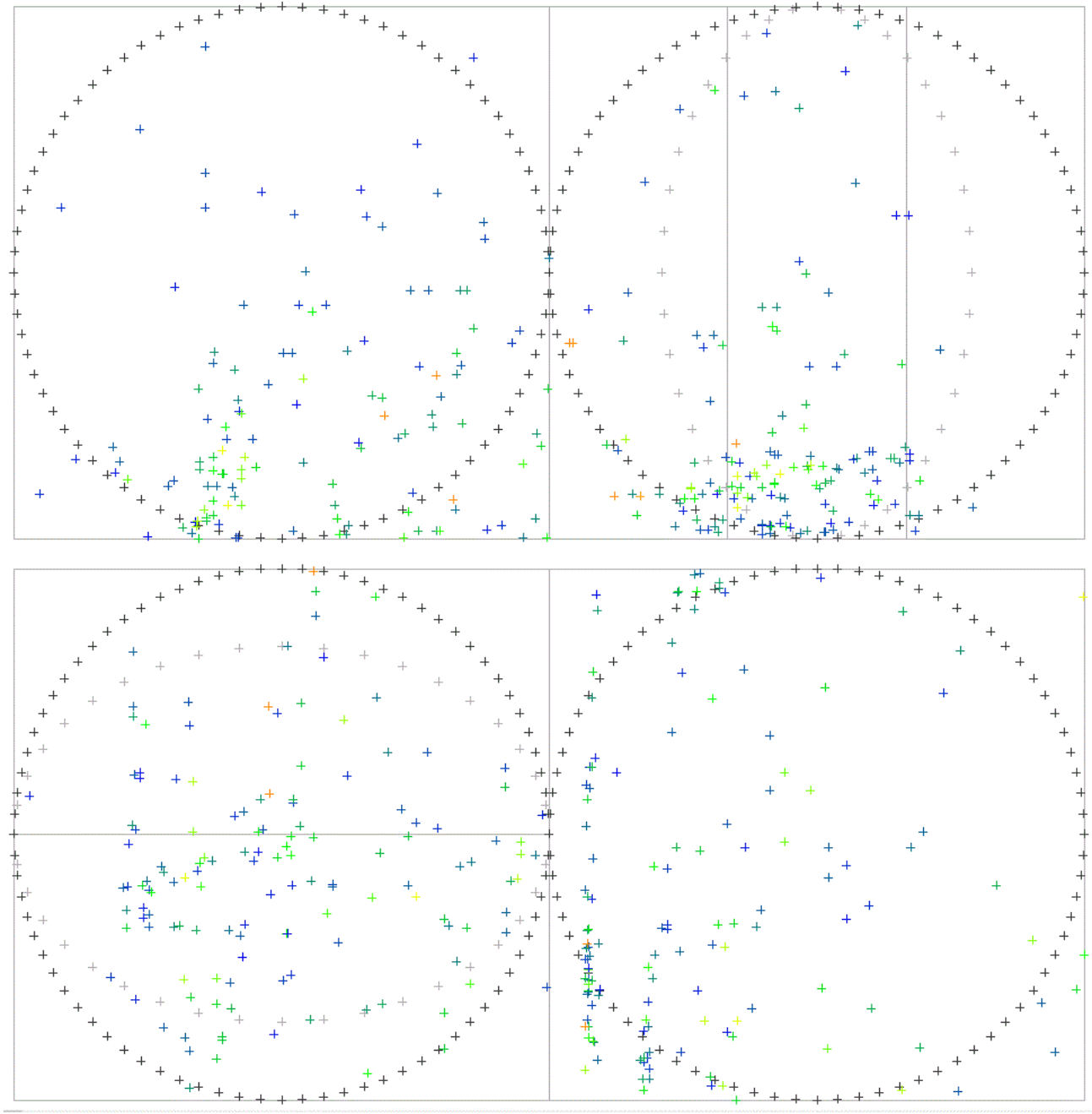

**Figure 14.** Points determined by the swarm along axes stage after a later, intermediate step, with *T*=0.25.

We followed instead a different scheme, in fact we perform line minimization along the same coordinate direction for all points of the swarm at each step! This means that we do something very similar to what we do in the first stage, the only difference is that instead of searching along a random direction in parameter space we use a coordinate direction (for computational efficiency, similarly to what we do in the first stage, line minimization may be skipped in case the first move does not improve over the starting point).

At first sight one may think that the two schemes are equivalent, but it is not so. The stacks that are generated are different, points in stack at each step are not simply substituted, they are "added", so that the worst point in stack will be eliminated at each step, therefore the differences after several steps can be great. Improving all points in stack at each step, rather than fully optimizing one point at a time, will produce more balanced stacks, so we preferred this scheme.

In any case performing line minimization along coordinate directions allows for some simplifications, the most notable one is that the line segment contained in the domain is always of length one, so that if $x_k$ is the current value of parameter $k$, the domain for *linmin* is quickly determined as *[-$x_k$, 1-$x_k$]*.

Fig. 14 shows the projections on the test planes of all points generated by this stage during the optimization process in one of the latest steps of the algorithm, when *T*=0.25.

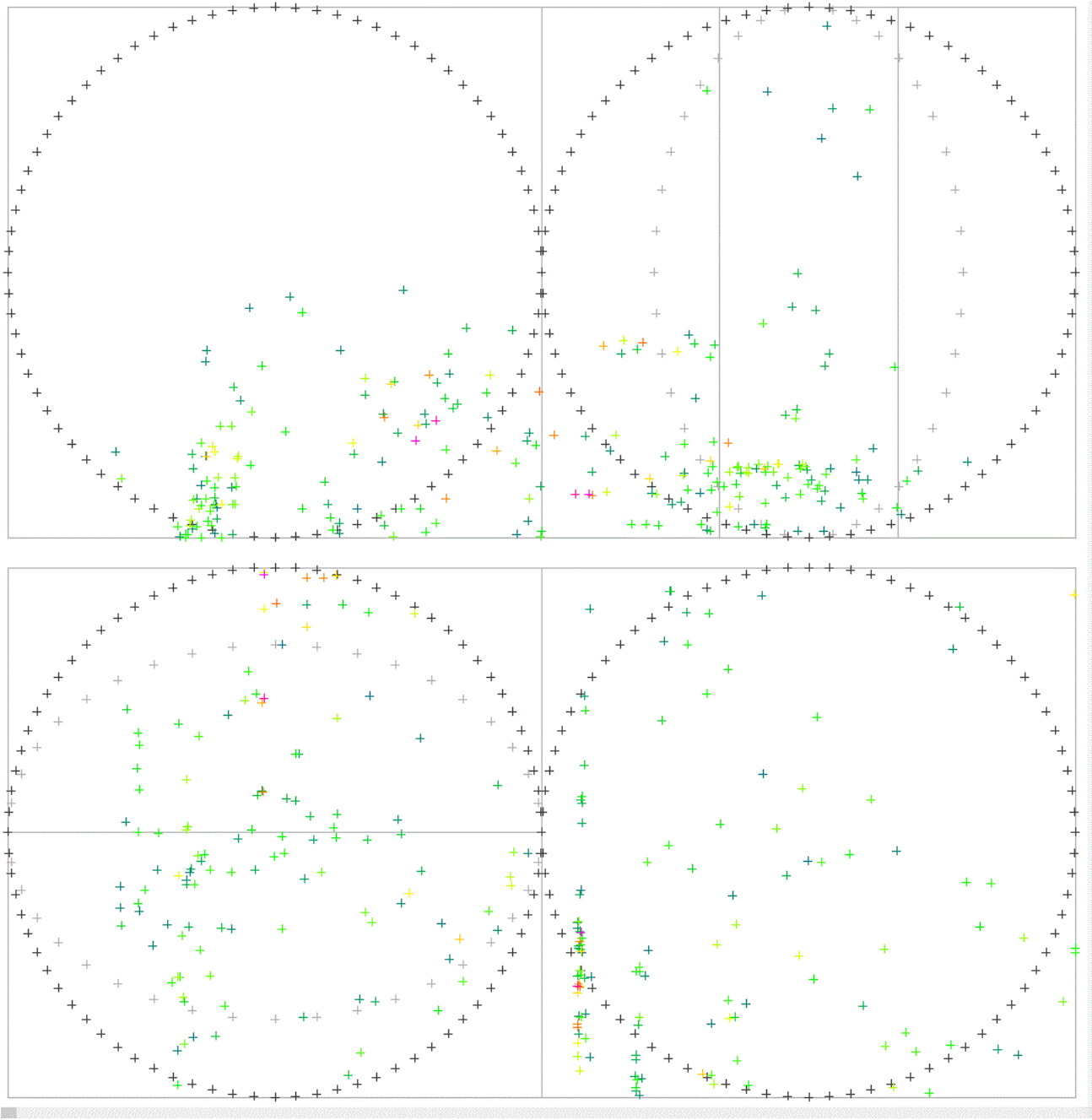

**Figure 15.** Projection over the test planes of the final points in stack, i.e. final result of the algorithm.

## 3    Final remarks

The final stack found by the algorithm is shown in Fig. 15. The best points found stand out in pseudocolors ranging from red to magenta, thus the best of them is recognized as the global optimum.

It is worth noting that for some applications not only the global optimum, but also the other sub-optimal solutions in stack may be of interest. For instance, when for different points the differences in the objective function are small, we may assume that these solutions are equivalent, therefore we may decide to choose any one of them randomly, or one in particular that we like, using in practice criteria which have not been embodied in the objective function.
Similarly, in case we are dealing with a true multimodal problem we should have multiple equivalent global optima. Our method, like other methods based on a swarm which gradually migrates towards the most important local minima, is naturally suited for this kind of problems [5], however, in case the number of optimal or near-optimal loci is large, it seems reasonable to employ also a large swarm to increase the probability of finding all of them.

## 4    Computational requirements

The computational requirements of the algorithm are obviously proportional to the complexity of the problem at hand, expressed by the complexity of the corresponding objective function. For regular simple functions of a few variables a few thousands of function evaluations can be enough for finding the optimal solutions, for complex and "wild" functions of several variables even hundred of thousands of function evaluations may be necessary. The test problems we dealt with can certainly be classified as complex, therefore the computational requirements are high. Typically we set the maximum number of function evaluations for each trial (composed of the four stages described above) at $10^4$, but this total is not evenly subdivided among the various stages, in fact about 30% more functions evaluations are assigned to the proximal search and swarm along axes stages, which are more demanding. A sequence (or set if executed concurrently) of ten trials is executed at the same temperature level, then this sequence is repeated at a lower temperature. Usually we perform five steps with $T$=1, 0.75, 0.5, 0.25, 0, therefore the total number of function evaluations at the end of execution of the algorithm sums up approximately to $5*10^5$. This number may seem too high, however we should remember that we dealt with problems in 11 dimensions, and that we computed the solutions with a precision (i.e. maximum error) on parameters of $10^{-4}$. Of course this precision can be achieved if sustained by meaningful variations of the function (because if the function is flat or almost flat any value in the neighborhood will do), moreover, because convergence in proximity of a minimum is usually quadratic, the final error can be much less than the required tolerance. To make a comparison, a dumb exhaustive search method of similar precision (with a maximum error of $(10^{-4})/2$ would require the computation of the function in $((10^4)^{11}) = 10^{44}$ points! The algorithm is thus very efficient with respect to exhaustive search.
Moreover, with respect to other methods which have been tested, imposing a similar computational burden, or even allowing for several times function evaluations, it should be pointed out that they did not find the global optimum solution most of the times, thus we claim the superiority of the proposed method.
The sequential execution of the algorithm (i.e. not parallelizing the execution of trials) in our test cases, required from about half an hour to a few hours of computation (depending on the cost of evaluating the function in the various cases) on an i7 machine.

## 5    Conclusions

We have presented an algorithm for the optimization of arbitrary functions in multidimensional space, which combines different strategies in a coherent scheme to take advantage of the different characteristics of each strategy. The algorithm is cascaded in four stages which are executed sequentially and repeated lowering "temperature" a prefixed number of times. Each stage embodies a different strategy, or otherwise said method. The methods which have been implemented are generally inspired by known techniques, which however have been substantially re-elaborated and improved, incorporating several original ideas, developed and tested in this work. The algorithm is suitable for attacking the most challenging optimization problems, e.g. problems characterized by a highly irregular objective function, presenting wild oscillations resulting in many local, sub-optimal points. It can manage multimodal problems as well.
In our experience it succeeded where other simpler methods failed. Every stage of the algorithm is substantially a stochastic process which depends on several sequences of (pseudo)random numbers following different suitable probability distributions (most of them proposed and described here for the first time), which therefore determine different random decisions and configurations of the solution space. To increase the probability of finding the global optimum, particularly in cases where this optimum in confined in a small region of the parameter space, the algorithm plans the execution of several "trials" (a trial is the execution of a complete sequence of all stages) with different random sequences. The execution of trials is naturally suitable for parallel execution.

**Acknowledgments**

The author wish to thank Tiziano Bombardi, for his suggestions after revising the draft of this document.

(*) Glauco Masotti was born in Ravenna, Italy, in 1955. He graduated summa cum laude Doctor in Electronic Engineering at the University of Bologna, Italy, in 1980.
Since then he has been involved in the design and development of various software systems, but he also got patents for some electronic devices.
He specialized in CAD/CAM systems, and in research in the field. Most of his work, being proprietary, has not been published. In early years, in Bologna, he cooperated at the development of the GBG drafting system for CAD.LAB (former name of Think3).
For COPIMAC he led the development of Aliseo, one of the first highly interactive, parametric and feature-based, solid modelers.
In 1989 he was visiting scholar at the University of Southern California (USC), in Los Angeles, where he worked at the design of an object-oriented geometric toolkit. He also proposed an extension to the C++ language, known as EC++ (Extended C++).
Back in Italy he undertook research at the Dept. of Electronics and Computer Science of the University of Bologna on symbolic methods and numerical problems in geometric algorithms. He also taught a course on Industrial Automation, with emphasis on geometric modeling and CAD/CAM integration.
In 1992 he joined Ecocad (later acquired by Think3), in Pesaro, Italy, where he was in charge of the design and development of various parts of Eureka, an advanced surface and solid modeler. One of his major contributions was the development of a module for assembly modeling and kinematic simulation, supporting the positioning of parts via mating constraints. In 1997 he was co-founder of IeS, Ingegneria e Software Srl, in Bologna, where he worked till 2001, developing software for automatic conversion of drawings in solid models, automatic production of exploded views, interpolation of measured 3D points, plane development of surfaces with minimal distortion, and finite element analysis. He was also involved in structural analysis and CFD problems.
In latest years, while trading for a living, he undertook autonomous researches in various fields, mainly on mathematical optimization and digital filters.
He is an enthusiast of sailing, which he practiced at a competitive level until 2002, racing on catamarans.